\input amstex
\documentstyle{amsppt}
\input bull-ppt
\keyedby{bull244/pah}

\define\ga{\alpha}
\define\gd{\delta}
\define\gy{\cyeps}
\define\gp{\pi}
\define\gu{\omega}
\define\gU{\Omega}
\define\fO{\scr O}
\define\fX{\scr X}
\define\fY{\scr Y}

\topmatter
\cvol{26}
\cvolyear{1992}
\cmonth{Jan}
\cyear{1992}
\cvolno{1}
\cpgs{113-118}
\ratitle
\title Lifting of cohomology and unobstructedness\\ 
of certain holomorphic maps\endtitle
\author Ziv Ran\endauthor
\shorttitle{Lifting and unobstructedness}
\address Institut des Hautes \'Etudes Scientifiques, 
Paris, France\endaddress
\curraddr Department of Mathematics, University of 
California,
Riverside, California 92521\endcurraddr
\date August 27, 1990\enddate
\subjclassrev{Primary\!\ 32G05,\!\ 32J27; Secondary 
32G20, 14J40}
\thanks Supported in part by NSF and IHES\endthanks
\abstract Let $f$ be a holomorphic mapping between
compact complex manifolds. We give a criterion for
$f$ to have {\it unobstructed deformations}, i.e.\ for
the local moduli space of $f$ to be smooth: this says,
roughly speaking, that the group of infinitesimal
deformations of $f$, when viewed as a functor, itself
satisfies a natural lifting property with respect to
infinitesimal deformations. This lifting property is
satisfied e.g.\ whenever the group in question admits a
`topological' or Hodge-theoretic interpretation, and we
give a number of examples, mainly involving Calabi-Yau
manifolds, where that is the case.\endabstract
\endtopmatter

\document
One of the most important objects associated to a compact
complex manifold $X$ is its {\it versal deformation} or
{\it Kuranishi family}
$$
\gp\:\fX\to\roman{Def}(X);$$
this is a holomorphic mapping onto a germ of an analytic
space $(\roman{Def}(X),0)$ (the Kuranishi space) with the
universal property that $\gp^{-1}(0)=X$ and that any
sufficiently small deformation of $X$ is induced by
pullback from $\gp$ by a map unique to 1st order. In
general, $\roman{Def}(X)$ is singular and even nonreduced;
in case $\roman{Def}(X)$ is smooth, i.e.\ a germ of the
origin in $\Bbb C^N$, we say that $X$ is {\it 
unobstructed}.
In an analogous fashion, a holomorphic mapping
$$
f\:X\to Y$$
also possesses a versal deformation, which in this case is
a diagram
$$
\matrix\format\c&\quad\c&\quad\l\\
\tilde f\:\fX & \longrightarrow & \quad\fY\\
\qquad\searrow & & \swarrow\\
& \roman{Def}(f)
\endmatrix$$
with a similar universal property. Again we say that $f$ is
unobstructed if $\roman{Def}(f)$ is smooth.

Now in \cite{R3}, we gave a criterion which deduces the
unobstructedness of a compact complex manifold $X$ from a
lifting property (in particular, deformation invariance)
of certain cohomology groups associated to $X$; this
implies in particular the unobstructedness of Calabi-Yau
manifolds, i.e.\ K\"ahler manifolds with trivial canonical
bundle $K_X$ (theorem of Bogomolov-Tian-Todorov
\cite{B, Ti, To}), as well as
that of certain manifolds with ``big'' anticanonical bundle
$-K_X$. In this note we announce an extension of our
criterion to the case of holomorphic maps of manifolds and
discuss some applications, mainly to maps whose source is
a Calabi-Yau manifold.

\heading 1. Generalities\endheading
Given a holomorphic map
$$
f\:X\to Y$$
of complex manifolds, we defined in \cite{R1} certain 
groups
$T^i_f$, $i\geq 0$, which are related to deformations of
$f$; in particular, $T^1_f$ is the group of 1st-order
deformations of $f$. For our present purposes, it will be
necessary to consider the corresponding relative groups
$T^i_{\tilde f/S}$, which are associated to a diagram
$$
\matrix\format\c&\quad\c&\quad\l\\
\tilde f\:\fX & \longrightarrow & \ \ \fY\\
\quad\searrow & & \swarrow\\
& S\endmatrix$$
with $\fX/S$, $\fY/S$ smooth (we call such a map $\tilde f$
an $S$-map, or a deformation of $f\,)$. In the notation of
\cite{R1, R2}, we have
$$
T^i_{\tilde f/S}=\roman{Ext}^i(\gd_1,\gd_0)$$
where $\gd_0\:f^*\fO_\fY\to\fO_\fX$, 
$\gd_1\:f^*\gU_{\fY/S}\to
\gU_{\fX/S}$ are the natural maps. As in \cite{R1}, we 
have an
exact sequence
$$
\aligned
0&\to T^0_{\tilde f/S}\to T^0_{\fX/S}\oplus T^0_{\fY/S}\to
\roman{Hom}_{\tilde f}(\gU_{\fY/S},\fO_\fX)\\
&\to T^1_{\tilde f/S}\to T^1_{\fX/S}\oplus T^1_{\fY/S}\to
\roman{Ext}^1_{\tilde f}(\gU_{\fY/S},\fO_\fX)\to\cdots
\endaligned\tag 1.1$$
where $T^i_{\fX/S}=H^i(T_{\fX/S})$, $T_{\fX/S}$ being the
relative tangent bundle and similarly for $T_{\fY/S}^i$,
$\roman{Hom}_{\tilde 
f}(\cdot,\cdot)=\roman{Hom}_\fX(\tilde f^*
\cdot,\cdot)$ and $\roman{Ext}^i_{\tilde f}(\cdot,\cdot)$
are its derived functors.

Now put $S_j=\roman{Spec}\,\Bbb C[\gy]/(\gy^j)$.
Our main general result, which is an analogue for
maps of a result given in \cite{R3} for manifolds,
is the following

\thm{Theorem-Construction 1.1} Suppose given $X_j/S_j$,
$Y_j/S_j$ smooth and $f_j\:X_j\to Y_j$ an $S_j$-map,
for some $j\geq 2$, and let $X_{j-1}/S_{j-1}$, $Y_{j-1}/
S_{j-1}$, $f_{j-1}\:X_{j-1}\to Y_{j-1}$ be their
respective restrictions via the natural inclusion $S_{j-1}
\hookrightarrow S_j$. Then

{\rm (i)} associated to $f_j$ is a canonical element
$\ga_{j-1}\in T^1_{f_{j-1}/S_{j-1}};$

{\rm (ii)} given any element $\ga_j\in T^1_{f_j/S_j}$ 
which maps to
$\ga_{j-1}$ under the natural restriction map 
$T^1_{f_j/S_j}\to
T^1_{f_{j-1}/S_{j-1}}$, there are canonically associated to
$\ga_j$ deformations $X_{j+1}/S_{j+1}$, $Y_{j+1}/S_{j+1}$ 
and
an $S_{j+1}$-map $f_{j+1}\:X_{j+1}\to Y_{j+1}$, 
extending $X_j/S_j$, $Y_j/S_j$ and $f_j\:X_j\to Y_j$
respectively.
\ethm

The proof is analogous to that of Theorem 1 in \cite{R3}
and will be presented elsewhere. In view of this theorem
it makes sense to give the following

\dfn{Definition 1.2} A map $f\:X\to Y$ is said to satisfy 
the
$T^1$-lifting property if for any deformation 
$f_j\:X_j/S_j\to
Y_j/S_j$ of $f$ and its restriction 
$f_{j-1}\!\:X_{j-1}/S_{j-1}\!
\to Y_{j-1}/S_{j-1}$, the natural map
$$
T^1_{f_j/S_j}\to T^1_{f_{j-1}/S_{j-1}}$$
is surjective.
\enddfn

Abusing terminology somewhat, we will say that $T^1_f$
is {\it deforma\-tion-invariant} if the groups 
$T^1_{f_j/S_j}$
are always free $S_j$-modules and their formation commutes
with base-change. Note, trivially, that whenever $T^1_f$
is deforma\-tion-invariant, $f$ satisfies the $T^1$-lifting
property. As an easy consequence of Theorem 1.1, we have 
the
following

\thm{Criterion 1.3} Suppose $f\:X\to Y$ is a map of
compact complex manifolds satisfying the $T^1$-lifting
property \RM(e.g.\ $T^1_f$ is deformation-invariant{\rm );}
then $f$ is unobstructed.
\ethm

\rem{Remark {\rm 1.4}} Various variants of this criterion
are possible, e.g.\ for deformations of maps $f\:X\to Y$
with fixed target $Y$. In the special case that $f$
is an embedding, with normal bundle $N$, we obtain that the
Hilbert scheme of submanifolds of $Y$ is smooth at the 
point
corresponding to $f(X)$ provided $H^0(N)$ satisfies the 
lifting property (e.g.\ is deformation-invariant).
Also, the {\it converse} to Criterion 1.3 is trivially 
true,
though we shall not need this.
\endrem

\heading 2. Applications\endheading
Unless otherwise specified, all spaces $X,Y$ considered 
here
are assumed smooth.

\thm{Theorem 2.1} Let $X$ be a Calabi-Yau manifold and
$f\:Y\hookrightarrow X$ the inclusion of a smooth
divisor. Then $f$ is unobstructed and moreover the image 
and
fibre of the natural map $\roman{Def}(f)\to\roman{Def}(X)$
are smooth.\ethm

\demo{Proof} In this case we may identify $T^1_f$ with
$H^1(T')$ where $T'$
is defined by the exact sequence
$$
0\to T'\to T_X\to N_{Y/X}\to 0,\tag 2.1$$
and it will suffice to prove deformation invariance of 
$H^1(T')$.
Now identifying $T_X\cong\gU_X^{n-1}$, $N_{Y/X}\cong\gU_Y
^{n-1}$, $n=\dim X$, we may write the cohomology sequence 
of
(2.1) as
$$
0\to H^{n-1,0}(Y)\to H^1(T')\to 
H^{n-1,1}(X)\overset{f^*}\to\to
H^{n-1,1}(Y)\cdots\.$$
As $H^{n-1,0}(Y)$ and $\roman{ker}(f^*)$ are both 
deformation-invariant,
so is $H^1(T')$, hence $f$ is unobstructed, and since 
moreover
the former groups are the respective tangent spaces to the
fibre and image of $\roman{Def}(f)\to\roman{Def}(X)$, the
latter are smooth.\QED\enddemo

A similar argument can be used to reprove a recent theorem
of C.~Voisin \cite{V} (see {\it op.\ cit.}\ for examples 
and
further results):

\thm{Theorem 2.2 {\rm (Voisin)}\rm} Let $X$ be a K\"ahler
symplectic manifold, with \RM(everywhere nondegenerate\RM)
symplectic form $\gu\in H^0(\gU^2_X)$, and $f\:Y\to X$
a Lagrangian embedding, i.e.\ $f^*\,\gu=0$ and $\dim Y=
\tfrac 12\dim X$. Then $f$ is unobstructed and the image 
and
fibre of the natural map $\roman{Def}(f)\to\roman{Def}(X)$
are smooth.
\ethm

\demo{Proof} In this case we may identify $T_X\cong\gU_X$,
$N_{Y/X}\cong\gU_Y$, and we may argue as in the proof of 
Theorem
2.1 (note that this property of being Lagrangian is
{\it open}).
\enddemo

Next we consider deformations of fibre spaces $f\:X^n\to 
Y^m$
with $X$ Calabi-Yau (i.e.\ $f$ is a flat map whose fibres
are reduced and connected). 
Note that for a fibre space $f$, its general fibre
is clearly a Calabi-Yau manifold. Also, it follows
easily from the sequence (1.1) that $\roman{Def}(f)
\hookrightarrow\roman{Def}(X)$. When $R^1f_*\fO_X=0$,
the morphism $\roman{Def}(f)\to\roman{Def}(X)$ is an
isomorphism by a theorem of Horikawa \cite{H}, hence
in that case unobstructedness of $f$ follows from that of
$X$. We will consider here two extreme cases: namely 
$m=n-1$ and
$m=1$.

\thm{Theorem 2.3} Let $f\:X\to Y$ be an elliptic
fibre space \RM(i.e.\ general fibre elliptic curve\RM)
with $X$ Calabi-Yau. Then $f$ is unobstructed.
\ethm

\demo{Proof} Using the usual exact sequence (1.1) and
Criterion 1.3, it suffices to prove the deformation 
invariance of
$$
\roman{ker}(H^1(T_X)\overset\ga\to\to 
H^0(Y,R^1f_*\fO_X\otimes T_Y))\.$$
Now by relative duality we have
$$
R^1f_*\fO_X\cong \gu^{-1}_{X/Y}\cong\gu_Y,$$
hence we may identify $\ga$ with the push-forward map (or
``integration over the fibre'')
$$
H^{n-1,1}(X)\to H^{n-2,0}(Y),$$
and in particular $\ker\ga$ is deformation-invariant.
(Note that we have $\roman{Def}(f)\cong\roman{Def}(X)$
whenever $\ga=0$, e.g.\ $H^{n-2,0}(Y)=0$, which holds
whenever $H^{n-2,0}(X)=0$.)
\enddemo

\thm{Theorem 2.4} Let $f\:X\to C$ be a fibre space
from a Calabi-Yau manifold to a smooth curve. Then
$f$ is unobstructed.
\ethm

\demo{Proof} Note that for any fibre $Y$ of $f$ we have
$$
h^0(\fO_Y(Y))=h^0(\fO_Y)=1,$$
and it follows that the scheme $\roman{Div}^0(X)$ 
parametrizing
reduced connected effective divisors of $X$ is smooth and
1-dimensional locally at the point corresponding to $Y$.
Consequently if we denote by
$$
p\:Z\to \roman{Div}^0(X)$$
the universal family and $q\:Z\to X$ the natural map, 
then we have
in fact a 1-1 correspondence between morphisms $f\:X\to C$
as above and smooth compact connected 1-dimensional 
components
$C\subset\roman{Div}^0(X)$ such that $q|p^{-1}(C)$ is an
isomorphism. Now it follows from Theorem 2.1 and its proof
that for any smooth fibre $Y$ of $f$, the locus 
$D'\subset \roman{Def}(X)$
of deformations over which $Y$ extends is smooth and
{\it independent} of $Y$. It follows that almost all, 
hence all, of
$C$ as component of $\roman{Div}^0(X)$ in fact extends over
$D'$, hence so does $f$, so that $D'=\roman{Def}(f)$, 
proving
the theorem.
\enddemo

In the intermediate cases, we have only much weaker 
results:

\thm{Theorem 2.5} Let $f\:X\to Y$ be a smooth morphism
and assume either

{\rm (i)} $K_X$ is trivial\RM; or

{\rm (ii)} $K_{X/Y}$ is trivial.

\noindent Then $\roman{Def}(f)\to\roman{Def}(Y)$ has 
smooth fibres.
\ethm

\demo{Proof} We will prove (ii), as (i) is similar. It 
suffices
to prove the deformation invariance of $H^1(T_{X/Y})$, 
where
$T_{X/Y}$ is the relative (vertical) tangent bundle. Now we
have
$$
T_{X/Y}\cong \gU^{n-1}_{X/Y}\otimes K^{-1}_{X/Y}\cong
\gU^{n-1}_{X/Y}\qquad n=\roman{dim}(X/Y)\.$$
By relative Hodge theory, $H^1(\gU^{n-1}_{X/Y})$ is a
direct summand of
$H^n(f^{-1}\fO_Y)$, and it will suffice to
prove the deformation invariance of the latter. We have a
Leray spectral sequence 
$$
H^p(Y,R^qf_*f^{-1}\fO_Y)\Rightarrow 
H^n(f^{-1}\fO_Y)\.\tag 2.2$$
However 
$H^p(Y\!,R^qf_*f^{-1}\fO_Y)=H^{p,0}(Y\!,R^qf_*\Bbb C_X)$
is a direct summand of $H^p(Y,R^qf_*\Bbb C_X)$,
hence the degeneration of the Leray spectral sequence of
$\Bbb C_X$ implies that of (2.2), hence the deformation
invariance of $H^n(f^{-1}\fO_Y)$.
\enddemo

\heading Acknowledgment\endheading
I am grateful to P.~Deligne for some helpful comments
concerning \cite{R3}, and to the IHES and Tel-Aviv
University, in particular Professor M.~Smorodinski, for
their hospitality.

\heading Added in proof\endheading
The above ideas are pursued further in the author's
preprints, {\it Hodge theory and the Hilbert scheme}
(September 1990) and {\it Hodge theory and deformations
of maps} (January 1991).


\Refs
\widestnumber\key{R1}
\ref\key B \by F. A. Bogomolov \paper Hamiltonian K\"ahler
manifolds \jour Dokl. Akad. Nauk SSSR \vol 243 \yr 1978 
\pages 1101--1104
\endref
\ref\key H \by E. Horikawa \paper Deformations of 
holomorphic
maps. {\rm III} \jour Math. Ann. \vol 222 \yr 1976 \pages 
275--282
\endref
\ref\key R1 \by Z. Ran \book Deformations of maps 
\bookinfo Algebraic
Curves and Projective Geometry (E. Ballico and C. 
Ciliberto, eds.),
Lecture Notes in Math. \publ vol. 1389, Springer-Verlag,
Berlin, 1989\endref
\ref\key R2 \bysame \paper Stability of certain 
holomorphic maps
\jour J. Differential Geom. \vol 34 \yr 1991 \pages 
37--47\endref
\ref\key R3 \bysame \paper Deformations of manifolds with 
torsion or
negative canonical bundle \jour J. Algebraic Geom. (to 
appear)\endref
\ref\key Ti \by G. Tian \paper Smoothness of the 
universal deformation
space of compact Calabi-Yau manifolds and its 
Peterson-Weil metric
\inbook Math. Aspects of String Theory (S. T. Yau, ed.) 
\pages 629--646,
World Scientific, Singapore, 1987\endref
\ref\key To \by A. N. Todorov \paper The Weil-Petersson 
geometry of the
moduli space of $SU(n\geq 3)$ \jour (Calabi-Yau) 
manifolds, preprint
IHES, November, 1988\endref
\ref\key V \by C. Voisin \paper Sur la stabilit\'e des 
sous-vari\'et\'es
Lagrangiennes des vari\'et\'es symplectiques holomorphes 
\jour
Orsay, preprint, April, 1990\endref


\endRefs
\enddocument